\documentclass[a4paper,12pt]{article}

\usepackage{amsmath}
\usepackage{amsthm}
\usepackage{amsfonts}

\parskip\smallskipamount
\newtheorem*{theorem*}{Theorem}
\newtheorem{theorem}{Theorem}

\newtheorem{remark}{Remark}

\newcommand{\pr}{\mathbb{P}}

\newcommand{\Z}{\mathbb{Z}}

\renewcommand{\S}{\mathcal{S}}

\renewcommand{\O}{\mathcal{O}}

\title{The randomly fluctuating hyperrectangles are spatially monotone  \footnote{
\noindent{Key-words: monotone occupied site probabilities; stochastic domination; symmetry; randomly fluctuating intervals.}} \vspace{-3mm}}   
\author{Achillefs Tzioufas\footnote{\small{\texttt{tzioufas@ime.usp.br}}}} 
\begin{document}

\maketitle
\vspace{-10mm}

\begin{abstract} 
We show that the probability of a site being occupied at any instance of time in the one-dimensional randomly fluctuating hyperrectangles processes decreases monotonically with respect to its distance from the origin. 
\end{abstract}
 
\section{Introduction} 

In this short note, we introduce a new class of Markovian spatial growth processes, which we call the randomly fluctuating hyperrectangles. We are interested in the property of the state of the process being at all times more likely to contain points that are nearer to the origin rather than points that are farther away. 
In a spatial stochastic process, a site is said to be occupied (or not) at a certain time according to whereas (or not) it is included at the state of the process at that time. Properties regarding occupied site probabilities are of inherent interest in the study of spatial stochastic processes, since these probabilities play a key r\^ole in their analysis and understanding, in addition to that it is natural to expect for this same reason their involvement in applications. Whereas occupied site probabilities are monotonically decreasing functions of their spatial coordinate at every fixed instance of time is an interesting in its own right problem that has been studied for various spatial stochastic processes, to which we refer to after mentioning our result below. 

We may briefly define the class of randomly fluctuating hyperrectangles processes 
informally as follows. The state-space of the process comprises of every hyperrectangle with vertices having integer-valued coordinates. The process evolves in discrete-time according to contraction phases  alternating with expansion phases. In the standard specifications case, contractions comprise of considering the sub-hyperrectangle obtained by sampling uniformly at random among all sub-hyperrectangles of the current hyperrectangle; whereas expansions comprise of shifting every face of the current hyperrectangle vertically in increasing direction according to independent geometrically distributed random values, and considering the sup-hyperrectangle formed by appropriately extending each shifted face. For instance, in two spatial dimensions an explicit definition of the process would be as follows. Let $\mathbf{R}$ be the set of (finite) rectangles with integer-coordinate vertices that, without loss of generality due to rotational invariance of the dynamics below, we assume to have sides parallel to the coordinate axes. Further, let $N(\zeta)$, $\zeta \in \mathbf{R}$, denote the coordinate value that points of the north side of $\zeta$ have in common, i.e.\ their projection on the vertical coordinate axis, and define $S(\zeta),E(\zeta),W(\zeta)$ analogously with regard to its south, east and west sides, respectively. Let also $\mathbf{R}(\zeta) = \{\xi \in \mathbf{R}:\xi \subseteq \zeta\}$, i.e.\ the $\sigma$-algebra of subsets of $\zeta$ in $\mathbf{R}$.
Let $(X^{i}(0), X^{i}(1), \dots)$, $i=N,S,E,W$ be independent collections of i.i.d.\ random variables such that $\pr(X^{i}(t)= n) = (1-p)p^{n-1}$, $n = 1,2, \dots$. The standard specifications randomly fluctuating rectangles process $(\zeta_{t}: t\geq 0)$ with parameter $p$ is a discrete-time Markov process on $\mathbf{R} \cup \emptyset$ the transition rates of which are determined as follows: given $\zeta_{t}$, we sample $\tilde{\zeta}_{t}$ uniformly at random\footnote{that is, according to the uniform probability measure, $\pr(\tilde{\zeta}_{t} = \zeta | \zeta_{t} = \eta) = \frac{1}{|\mathbf{R}(\eta)+1|}$, for all $  \zeta \in \mathbf{R}(\eta) \cup \emptyset$.} from $\mathbf{R}(\zeta_{t}) \cup \emptyset$ and, whenever $\tilde{\zeta}_{t} \not= \emptyset$, let $\zeta_{t+1}$ be such that $i(\zeta_{t+1}) = i(\tilde{\zeta}_{t}) + X^{i}_{t}$ for $i =N,E$, and that $i(\zeta_{t+1}) = i(\tilde{\zeta}_{t}) - X^{i}_{t}$ for $i= S,W$, otherwise let $\zeta_{k} = \emptyset$ for all $k \geq t+1$. Thus, contractions of this type comprise of sampling a sub-rectangle according to the probability measure which assigns equal mass to all elements, including the empty set; whereas expansions will always comprise of independent geometrically distributed outward-shifting of each side, and considering the sup-rectangle formed by joining them.

We may briefly summarize some of the intrinsic features exhibited by the randomly fluctuating hyperrectangles processes as follows. We note that the process in the standard specifications case lacks an obvious graphical representation, i.e.\ a coupling construction of versions of the process with different starting states, and hence do not fall into the general framework of interacting particle systems, or percolation processes. Further, regarding the restriction to finite initial states, we note that the uniform contraction rule is not otherwise well-defined, as an aftereffect of the mere fact that uniform (i.e.\  assigning equal probability to every element) probability measures on countable spaces are not compatible with the standard axioms of probability theory. We also note that it is not difficult to see that this class of Markov processes is irreducible, in that, for any state $\zeta$, there is $t$ such that $\pr(\zeta_{t}= \zeta)>0$.

Our Theorem \ref{main} stated in the next section regards the randomly fluctuating hyperrectangles class of  processes in one (spatial) dimension started from the single site at the origin. Although we work out the details in our proofs in the one-dimensional case only, we nevertheless find that our result and arguments extend analogously in any dimension.  Further,  whereas the specification of a geometric distribution for the expansion phases is essential here, the result and technique of proof apply directly for the randomly fluctuating hyperrectangles processes with other, qualitatively different  than uniform types of contraction specifications (see Remarks \ref{rem1B} and \ref{rem1C} below). We show in Theorem \ref{main} that at any instance of time the probability that a point is occupied decreases monotonically with respect to its distance from the origin.  We furthermore show in Theorem \ref{main}  that, at any fixed instance of  time, the occupied site probabilities is an even function; that is, that the probability that a site is occupied is equal to the probability that its symmetric about the origin counterpart site is occupied. 

We refer to studies regarding the corresponding spatial monotonicity property for other stochastic spatial growth processes as follows. In regard to the basic one-dimensional contact process, which is the continuous-time analog of two-dimensional oriented percolation, Gray \cite{G91} introduced the spatial monotonicity of its occupied site probabilities, among other intriguing properties regarding them. The detailed and elaborate proof of this notable result is given by Andjel and Gray \cite{AG14} and by Andjel and Sued \cite{AS08}.  Whereas the corresponding property holds for undirected percolation on integer lattices is in general an open question. A partial result in the direction of a positive reply is obtained by de Lima {\it et.al.} \cite{LPS15}. Whereas occupied site probabilities are monotone for general one-dimensional attractive spin-systems is also in the case of finite initial configurations an important open problem. Regarding first-passage percolation, Hammersley and Welsh \cite{HW65} raised the corresponding question of spatial monotonicity for first-passage times,  which also remains to date an open question. For a partial result in the direction of a negative reply, see van den Berg \cite{B83}; for partial results in the direction of a positive reply, see the more recent work by Gou\' er\' e \cite{G14} and the references therein. We also note that the corresponding stochastic monotonicity result regarding symmetric branching random walks is derived in Lemma 11 by Lalley and Zheng \cite{LZ11}, the proof of which relies crucially on the independence of the descendancy of distinct particles, due to permitting in this process for an arbitrary number of particles per site. 

The remainder of this note is organized as follows. We state our Theorem \ref{main} in Section \ref{statements} next. We give its proof in the subsequent Section \ref{Secproofs}.

\section{Statement of Theorem \ref{main}}\label{statements}

The {randomly fluctuating intervals} process is defined as follows. Let $\mathbf{I}(\zeta)$ be the set of all integer interval subsets of $\zeta \subseteq \Z$ including the empty set, where $\Z$ denotes the integers, and simply write $\mathbf{I}$ for $\mathbf{I}(\Z)$. Let further $R(\zeta) = \sup\zeta$ and $L(\zeta) = \inf \zeta$, with the convention that $\sup \emptyset=-\infty$. Furthermore, throughout here,  $(N^{L}_{t}: t \geq 0)$ and $(N^{R}_{t}: t \geq 0)$  will denote independent collections of i.i.d.\ geometric r.v.\ such that $\pr(N^{L}_{t}= n) = (1-p)p^{n-1}$, $n =1,2, \dots$, where $p \in (0,1)$ is called the {expansion parameter}.
The {standard specifications randomly fluctuating intervals} is a Markov process  $(\zeta_{t}: t\geq 0)$ on $\mathbf{I}$ with parameter $p \in (0,1)$ defined as follows: Given $\zeta_{t}$, choose $\tilde{\zeta}_{t}$ uniformly at random from $\mathbf{I}(\zeta_{t})$ and, whenever $\tilde{\zeta}_{t} \not= \emptyset$, set $\zeta_{t+1} = \{L(\tilde{\zeta}_{t}) - N^{L}_{t}, \dots , R(\tilde{\zeta}_{t}) + N^{R}_{t} \}$, otherwise set $\zeta_{k} = \emptyset$ for all $k \geq t+1$. To state our theorem next, let $\zeta_{t}^{O}$ be the standard specifications randomly fluctuating intervals process started at the origin, and let also $f_{t}(x) = \pr(x \in \zeta_{t}^{O})$. 

\begin{theorem}\label{main}
\[
\mbox{For all $t$, $f_{t}(x)$ is an even function that is decreasing in $|x|$}. 
\]
\end{theorem}

\begin{remark}\label{rem1B}  
\normalfont The arguments in the proof can be adapted to also establish this result for {the randomly fluctuating intervals with general contraction}. We describe this Markov process $(\eta_{t})$ on $\mathbf{I}$. Given  $\eta_{t}$, sample $\tilde{\eta}_{t}$ uniformly at random among intervals in $\mathbf{I}(\eta_{t})$ with size $X \sim \phi(k;|\eta_{t}|)$, where $\phi(k;n)$ is an arbitrary probability mass function of a discrete random variable assuming values $k = 0,\dots,n$ and,  whenever $\tilde{\eta}_{t} \not= \emptyset$, set $\eta_{t+1} = \{L(\tilde{\eta}_{t}) - N^{L}_{t}, \dots , R(\tilde{\eta}_{t}) + N^{R}_{t} \}$, otherwise set $\eta_{k} = \emptyset$ for all $k \geq t+1$. 
\end{remark}

\begin{remark}\label{rem1C}
\normalfont Regarding other variations, we may also allow for transitions of $(\zeta_{t})$ to the empty set at time $t$ to occur with probability $p_{\emptyset} = p_{\emptyset}(p, |\zeta_{t}|)$ and sampling $\tilde{\zeta}_{t}$ uniformly at random from $\mathbf{I}(\zeta_{t}) \backslash \{\emptyset\}$ with probability $1-p_{\emptyset}$. Additionally, we may substitute the uniformly at random rule by sampling instead
uniformly at random with repetition from $\{x:x \in \zeta_{t}\}$ the endpoints of $\tilde{\zeta}_{t}$.  
\end{remark}

The technique of proof of the part of Theorem \ref{main} that regards the occupancy function being even relies on a simple coupling construction of two doppelg\"anger (mirror images about the origin) versions of the processes. This part of the result is in fact an elementary consequence of the symmetry inherent in the definition of the process, and we include the proof for completeness. That of the second part relies on a different and much more elaborate coupling of this type that allows us to control the spatial competition of two processes about the origin until the so called {coupling time}, which may well be infinite. We believe that directly analogous extensions of the arguments involved can be used to yield the result in higher dimensions (in $L^{1}$ distance).


\section{Proof of Theorem \ref{main}}\label{Secproofs}

\begin{proof}[Proof of Theorem \ref{main}]

The statement is equivalent to showing that 
\begin{equation}\label{one}
\pr(x \in \zeta_{t}^{O}) \geq \pr(x+1 \in \zeta_{t}^{O}), \mbox{ for all } x \geq 0,  
\end{equation}
and that 
\begin{equation}\label{one2}
\pr(x \in \zeta_{t}^{O}) = \pr(-x \in \zeta_{t}^{O}).
\end{equation}

We first show (\ref{one2}). To do this we will exploit  innate symmetries of the definition of the process to construct on the same probability space two copies of it, $(\zeta_{t}^{O})$ and $(\eta_{t}^{O})$, such that 
\begin{equation}\label{two}
\zeta_{s}^{O} = - \eta_{s}^{O}
\end{equation}
for all  $s$, where, given integer interval $i = \{a, \dots, b\}$, $b \geq a$, $i \not= \emptyset$, we let $-i$ denote its {reflection} about the origin integer interval,  $-i := \{-b, \dots, -a\}$,  where by convention $-\emptyset:= \emptyset$.  Note that this conclusion implies in particular that $\{ x \in \zeta_{s}^{O} \}$ if and only if $\{ - x \in \eta_{s}^{O} \}$, and hence it implies $(\ref{one2})$.  Since  $(\ref{two})$  is obviously true when $s=0$, we assume that it holds for some $s=t$ and prescribe coupled transitions that imply it for $s=t+1$. In case $\zeta_{t}^{O} = - \eta_{t}^{O} = \emptyset$ there is nothing to prove as the empty set is absorbing, and hence we may assume  $\zeta_{t}^{O} = - \eta_{t}^{O} \not= \emptyset$. Note that we may choose contractions for $\eta_{t}^{O}$ to be reflections about the origin of those for $\zeta_{t}^{O}$ with the correct marginal transition rates since also the cardinal of $\zeta_{t}^{O}$ and of $-\eta_{t}^{O}$ are equal; that is, given the contraction $\tilde{\zeta}_{t}^{O}$ for $\zeta_{t}^{O}$, we choose the contraction for $\eta_{t}^{O}$, $\tilde{\eta}_{t}^{O}$, such that $\tilde{\zeta}_{t}^{O} = -\tilde{\eta}_{t}^{O}$. The prescription of the expansion phases in this coupling is such that if expansions for $\zeta_{t}^{O}$ are defined by means of $N^{L}_{t}$ and $N^{R}_{t}$, those of $\eta_{t}^{O}$ are obtained by interchanging their roles, that is, we use $N^{L}_{t}$ and $N^{R}_{t}$ for the right-end and left-end expansions respectively. This  implies $(\ref{two})$ for $s=t+1$, hence proving  $(\ref{one2})$.

We now turn to the proof of $(\ref{one})$. Let $\zeta_{t}^{-1}$ denote the process with $\zeta_{t=0}^{-1} = \{-1\}$.  Due to that $(\zeta_{t}^{-1})$ is a version of $(\zeta_{t}^{O})$ shifted to the left by one site, it is easily seen that $\pr(x \in \zeta_{t}^{-1}) = \pr(x+1 \in \zeta_{t}^{O})$, so that it suffices to show that
\begin{equation}\label{des0}
\pr(x \in \zeta_{t}^{O}) \geq \pr(x \in \zeta_{t}^{-1}),  \mbox{ for all } x \geq 0. 
\end{equation}

To carry out this coupling construction, we will evoke below the following equivalent surface description of the expansion phase. Let $(\omega_{t}^{L}: t \geq 1)$  be a collection of collections of i.i.d.\  Bernoulli r.v.\ $\omega_{t}^{L} := (\omega_{t}^{L}(n): n \geq 1)$ such that $\omega_{t}^{L}(n)=1$ with probability $p$, and let further $(\omega_{t}^{R}: t \geq 1)$ be an independent copy of  $(\omega_{t}^{L}: t \geq 1)$. Note that we may now equivalently define the process $\zeta_{t}^{O}$ with expansions according to $N_{t}^{R} := \min\{n: \omega_{t}^{R}(n) = 0\}$ and  $N_{t}^{L} := \min\{n: \omega_{t}^{L}(n) = 0\}$, and contractions as before.

Further, it will be seen below that it is convenient for what follows to relabel site $0$ as $1$, $1$ as $2$, etc., whereas sites $-1$, $-2$, etc.\ maintain their labeling, and hence site $0$ is omitted in this labeling now. Accordingly, we will also denote $\zeta_{t}^{O}$ simply by $\zeta_{t}^{+}$ and also denote $\zeta_{t}^{-1}$ simply by $\zeta_{t}^{-}$ below. Note that $\zeta_{t=0}^{+}=\{1\}$ and that  $\zeta_{t=0}^{-}=\{-1\}$ in the new labeling. We need to show that 
\begin{equation}\label{des1}
\pr(x \in \zeta_{t}^{+}) \geq \pr(x \in \zeta_{t}^{-}),  \mbox{ for all } x \geq 1,
\end{equation}
which is equivalent to (\ref{des0}) with the new labeling. 

We will also need the following definitions. Let $i$ and $j$ be integer intervals with respect to the new labeling, that is, sets of consecutive (w.r.t.\ the new labeling) integers. We say that $i$ and $j$ are {antithetic} if and only if $L(j) = -R(i)$ and $R(j) = -L(i)$, where also by convention the empty set is antithetic to itself. Further, let $\mathbf{A}$ be the collection of ordered pairs of integer intervals $(i, i')$ such that $i$ and $i'$ are antithetic {and} $i$ is such that $L(i) = -1, -2, \dots$ and that $R(i) \leq -L(i)$, in words, $i$ contains at least as many negative integers as positive ones and at least one negative one. In addition, we let $\mathbf{O}$ denote the subset of $\mathbf{A}$ containing the single pair $(\emptyset, \emptyset)$, and also let $\mathbf{S}$ denote the subset of $\mathbf{A}$ containing all $(i, i')$ such that $i = i'$ and $i \not=\emptyset$.

We shall prescribe a coupling for the two processes such that for any given $(\zeta_{t}^{-}, \zeta_{t}^{+}) \in \mathbf{A} \backslash \mathbf{S} \cup \mathbf{O}$ one of the following three mutually exclusive possibilities for the coupled transitions occurs. Either $(\zeta_{t+1}^{-}, \zeta_{t+1}^{+}) \in \O$, or  $(\zeta_{t+1}^{-}, \zeta_{t+1}^{+}) \in \S$, or $(\zeta_{t+1}^{-}, \zeta_{t+1}^{+}) \in \mathbf{A} \backslash \S \cup \O$. We note that the possibility of $(\zeta_{t+1}^{-}, \zeta_{t+1}^{+}) \in \S$ will arise both according to the expansion and the contraction coupled phases. Further, if $(\zeta_{t+1}^{-}, \zeta_{t+1}^{+}) \in \S$, then the coupling is such that $(\zeta_{k}^{-}, \zeta_{k}^{+}) \in \S \cup \O$ for all $k \geq t+1$. Hence, since $(\zeta_{t = 0}^{-}, \zeta_{t = 0}^{+}) \in \mathbf{A}  \backslash \mathbf{S} \cup \mathbf{O}$, this gives that we have constructed a coupling of $\zeta_{t}^{-}$ and $\zeta_{t}^{+}$ with the property that 
\begin{equation*}\label{des} 
\zeta_{t}^{+} \cap  \{1, 2, \dots \} \supseteq \zeta_{t}^{-} \cap  \{1, 2, \dots \}, \mbox{ for all } t,  
\end{equation*}
a.s., which implies (\ref{des1}), and hence showing such a coupling completes the proof.

We first prescribe the coupled contraction phases.  Given integer interval $i$, we denote the (unique) antithetic integer interval by $i^{T}$. Given integer interval $i$, let $\mathbf{F}(i)$ denote the set of all of its integer sub-intervals (including the empty set). 
Consider the one-to-one correspondence defined by the (bijective) function $\psi: \mathbf{F}(i) \rightarrow \mathbf{F}(i')$ given by mapping each interval to its antithetic one, and the empty set to the empty set, that is, the function $\psi$ such that $\psi(j) = j$ for all $j \subseteq i \cap i'$, and that $\psi(\emptyset) = \emptyset$, and further that $\psi(j) = j^{T}$, for all other $j$. The key observation here is that for any given $(i, i') \in \mathbf{A}$ and all $j \in \mathbf{F}(i)$, we have that $\psi(j)$ is such that either $(j, \psi(j)) \in \mathbf{A}$, or $j = \psi(j)$. 

Let $Z_{t} = \zeta_{t}^{-} \cap \zeta_{t}^{+}$, where possibly $Z_{t} = \emptyset$. The coupling of the contraction phases proceeds by sampling $\tilde{\zeta}_{t}^{-}$ and setting $\tilde{\zeta}_{t}^{+} = \psi(\tilde{\zeta}_{t}^{-})$. The fact that this coupling gives the correct marginals for the two processes is an immediate consequence of that $|\mathbf{F}(\zeta_{t}^{-})| = |\mathbf{F}(\zeta_{t}^{+})|$ and that $\psi$ defines a bijection among all possible contraction states. Note now that by the definition of $\psi$ we have the following three mutually exclusive possibilities. If $(\tilde{\zeta}_{t}^{-}, \tilde{\zeta}_{t}^{+}) \in \O$ then $(\zeta_{t+1}^{-}, \zeta_{t+1}^{+}) \in \O$ holds. If $(\tilde{\zeta}_{t}^{-}, \tilde{\zeta}_{t}^{+}) \notin \O$  then we either sample $\tilde{\zeta}_{t}^{-} \not\subseteq Z_{t}$, and in this case we have that $(\tilde{\zeta}_{t}^{-}, \tilde{\zeta}_{t}^{+}) \in \mathbf{A} \backslash \O \cup \S$, or we sample $\tilde{\zeta}_{t}^{-} \subseteq Z_{t}$, and in this case we have that $\tilde{\zeta}_{t}^{+} = \tilde{\zeta}_{t}^{-}$. 

In the case $\tilde{\zeta}_{t}^{+} = \tilde{\zeta}_{t}^{-}$ we proceed coupling transitions for the two processes identically for all future times, so that $(\zeta_{t+1}^{-}, \zeta_{t+1}^{+}) \in \S$ and further  $(\zeta_{k}^{-}, \zeta_{k}^{+}) \in \S \cup \O$ for all $k \geq t+2$. Hence the only remaining case for which we need to prescribe the expansion phases is that of $(\tilde{\zeta}_{t}^{-}, \tilde{\zeta}_{t}^{+}) \in \mathbf{A} \backslash \O \cup \S$. This prescription will be such that: either $(\zeta_{t+1}^{-}, \zeta_{t+1}^{+}) \in \S$  or $(\zeta_{t+1}^{-}, \zeta_{t+1}^{+}) \in \mathbf{A} \backslash \S \cup \O$. Again, in the former case, we proceed coupling transitions for the two processes identically for all future times. Thus, note that, with such prescription in hand, we have achieved the coupling of the two processes with all claimed properties.  

To carry out the final part of the construction for obtaining $(\zeta_{t+1}^{-}, \zeta_{t+1}^{+}$) in the case $(\tilde{\zeta}_{t}^{-}, \tilde{\zeta}_{t}^{+}) \in \mathbf{A} \backslash \O \cup \S$ we will exploit the equivalent surface description of the expansion phase by means of Bernoulli r.v.\ that was detailed in the paragraph following (\ref{des0}) above. That is, we will describe the two expansion phases by means of coupled collections of Bernoulli r.v.  Let $(\omega_{t+1}^{-, R}(n): n \geq1)$ and $(\omega_{t+1}^{-, L}(n): n \geq1)$ be associated with the right-side and left-side respectively expansion for  $\tilde{\zeta}_{t}^{-}$, and similarly $(\omega_{t+1}^{+, R}(n): n \geq1)$, $(\omega_{t+1}^{+, L}(n): n \geq1)$ will be associated with right-side and left-side respectively expansion for $\tilde{\zeta}_{t}^{+}$.

Let $g_{t} = R(\tilde{\zeta}_{t}^{+}) - R(\tilde{\zeta}_{t}^{-})$. We sample independent $\omega_{t+1}^{-, R}, \omega_{t+1}^{-, L}$ and consider the following two cases.  If the following condition is satisfied
\begin{equation}\label{eqg}
\omega_{t+1}^{-, R}(n) =1, \mbox{ for all } i = 1, \dots, g_{t}, 
\end{equation}
then we set $\omega_{t+1}^{+, L}(n) =  \omega_{t+1}^{-, R}(n)$, for all $n = 1, \dots, g_{t}$, and $\omega^{-, R}_{t+1}(g_{t}+n) =  \omega^{+, R}_{t+1}(n)$, as well as $\omega^{-, L}_{t+1}(n) = \omega^{+, L}_{t+1}(g_{t} + n)$, $n \geq 1$. Otherwise, we set $\omega_{t+1}^{-, R}(n) = \omega_{t+1}^{+, L}(n)$  and 
$\omega_{t+1}^{-, L}(n) = \omega_{t+1}^{+, R}(n)$, for all $n$.  Hence, in case (\ref{eqg}) checks we have that $\zeta_{t+1}^{+} = \zeta_{t+1}^{-}$, whereas in case not, we have that $(\zeta_{t+1}^{-}, \zeta_{t+1}^{+}) \in \mathbf{A} \backslash \S \cup \O$, as required. 
\end{proof}

\noindent {\it Acknowledgments}. Research for this work was carried through under the auspices of Conicet, Argentina. The author is currently supported by FAPESP fellowship (grant 2016/03988-5), part of FAPESP project Research, Innovation and Dissemination Center for Neuromathematics (grant 2013/ 07699-0).

\textsc{ \\
\noindent Instituto de Matem\' atica e Estat\' istica,\\
Universidade de S\~ ao Paulo\\
Rua do Mat\~ ao, 1010\\
CEP 05508-900- S\~ ao Paulo\\
Brasil}
\end{document}